\documentclass[
 reprint,  
 superscriptaddress,
 amsmath,amssymb,aps,showpacs,prl
]{revtex4-1}
\usepackage{graphicx}
\usepackage{dcolumn}
\usepackage{bm}
\usepackage[colorlinks=true,          
            linkcolor=blue,
            citecolor=blue,
            urlcolor=blue,
            ps2pdf,                   
            dvips,
            pdfpagelayout=OneColumn,  
            pdfstartview=FitH,        
            pdfstartpage=1,           
            pdfdisplaydoctitle=true,  
            pdfcreator  = {LaTeX with hyperref package},
            pdfproducer = {dvips + ps2pdf}
            ]{hyperref}
\hypersetup{
            pdfauthor   = {Matthias Heymann, Stephen Teitsworth, Jonathan Mattingly},
            pdftitle    = {Rare Transition Events in Nonequilibrium Systems with State-Dependent Noise: Application to Stochastic Current Switching in Semiconductor Superlattices},
            pdfsubject  = {PACS numbers: 05.40.-a, 05.10.-a, 05.10.Gg, 72.70.+m},
            pdfkeywords = {superlattice, current switching, current branches, shot noise, transition path, bistability, large deviation theory, geometric minimum action method, gMAM}
}
\let\url\nolinkurl 

\usepackage{dsfont}

\allowdisplaybreaks
\newcommand{\Vth}{V_{\mathrm{\it \!th}}}

\newcommand{\ND}{N_{\!D}}
\newcommand{\vM}{v_{\hspace{-.5pt}M}}

\newcommand{\R}{\mathbb{R}}
\newcommand{\Rn}{\R^N}

\newcommand{\Gx}{\Gamma_{x_1}^{x_2}}
\newcommand{\skp}[2]{\langle#1,#2\rangle}
\newcommand{\skpA}[2]{\langle#1,#2\rangle_x}
\newcommand{\One}{\mathds{1}}

\begin{document}
\preprint{APS/123-QED}


\title{\texorpdfstring{Rare Transition Events in Nonequilibrium Systems with State-Dependent Noise:\\Application to Stochastic Current Switching in Semiconductor Superlattices}{Rare Transition Events in Nonequilibrium Systems with State-Dependent Noise: Application to Stochastic Current Switching in Semiconductor Superlattices}}

\author{Matthias Heymann}
 \email{heymann@math.duke.edu}
 \homepage{www.matthiasheymann.de}
\affiliation{
  Duke University Mathematics Department, Durham, North Carolina 27708, USA
}

\author{Stephen W. Teitsworth}
 \email{teitso@phy.duke.edu}
\affiliation{Duke University Physics Department, Durham, North Carolina 27708, USA}

\author{Jonathan C. Mattingly}
 \email{jonm@math.duke.edu}
\affiliation{
  Duke University Mathematics Department, Durham, North Carolina 27708, USA
}

\date{\today}

\begin{abstract}
Using recent mathematical advances, a geometric approach to rare noise-driven transition events in nonequilibrium systems is given, and an algorithm for computing the maximum likelihood transition curve is generalized to the case of \textit{state-dependent} noise.
It is applied to a model of electronic transport in semiconductor superlattices to investigate transitions between metastable electric field distributions. When the applied voltage $V$ is varied near a saddle-node bifurcation at $\Vth$, the mean life time $\langle T\rangle$ of the initial metastable state is shown to scale like $\log\langle T\rangle\propto|\Vth-V|^{3/2}$ as $V\nearrow\Vth$.
\end{abstract}


\maketitle



The study of noise-driven transitions between metastable states is of current interest for a range of systems, for example, chemical reaction systems \cite{denOtter}, nano- and micromechanical oscillators \cite{AldridgeCleland,ChanStambaugh,ChanDykman}, magnetic tunnel junctions~\cite{Cui}, and biochemical networks~\cite{AllenWarren}. Accordingly, the development of mathematical and numerical methods \cite{String1,MAM,CPAM_Hey1,PRL_Hey1} for understanding in particular the mean transition (or escape) time and the most probable escape path is an active area with implications for a variety of fields. Since the widely-used string method \cite{String1} is restricted to the case of gradient drift, the geometric minimum action method (gMAM) \cite{CPAM_Hey1,PRL_Hey1,JCP_Hey1}
was recently introduced for non-gradient systems (as a successor of the minimum action method, MAM \cite{MAM}). However, none of its implementations can at present deal with multiplicative (i.e., state-dependent) noise, leaving it inapplicable to a number of important physical problems.

In this paper we fill this gap by extending gMAM to systems with multiplicative noise, and we also introduce a useful new method for locating unstable equilibrium points in high-dimensional systems.
The mathematical framework, large deviation theory \cite{WentzellFreidlin}, is presented in a non-standard geometric way (i.e., based on \textit{un}parameterized curves), which allows us to incorporate a key result of a recent study \cite{Existence1} proving the existence of a maximum likelihood transition curve~$\gamma^\star$.

These techniques are illustrated by applying them to the problem of stochastic current switching in a semiconductor superlattice (SL), which is a linear array of $N>1$ quantum wells that are coupled to each other via electron tunneling.
While the \textit{deterministic} nonlinear properties of SLs have been studied extensively and are generally well-understood \cite{BonillaTeitsworth}, very little theoretical work \cite{BonillaSanchezSoler} exists for explaining experimental measurements \cite{Rogozia} of noise-induced current switching in this system.
One object of interest is the mean escape time $\langle T\rangle$ from the metastable states in the SL system, and its dependence on key control parameters such as the applied voltage $V$.
Since the noise in the SL model is multiplicative, previous studies could not yet use gMAM and instead relied on direct simulation of the underlying stochastic differential equations (SDEs) \cite{BonillaSanchezSoler}. The resulting lack of accuracy and high computational cost (note that $\langle T\rangle$ becomes exponentially large in the zero-noise limit) can now be circumvented by the non-trivial extension of gMAM introduced here.

This allows us, for the first time, to determine the scaling behavior of $\langle T\rangle$ for the SL system.  Specifically, we show that $\log\langle T\rangle \propto |\Vth-V|^{3/2}$ as $V\nearrow\Vth$, where $\Vth$ is the value above which the initial metastable state ceases to exist. Previously, this scaling result was only established for the one-dimensional case, i.e., for the case of a single quantum well with small cross-sectional size~\cite{Tretiakov}.



\vspace{.2cm}
We begin by presenting the SDE model for the SL system.
Electronic transport properties of weakly-coupled superlattices are accurately described by a spatially discrete model that incorporates sequential resonant tunneling between successive quantum wells \cite{BonillaTeitsworth,XuTeitsworth}.
For a superlattice with $N$ quantum wells, this model consists of the Poisson and the charge continuity equations
\begin{align}
F_i-F_{i-1} &= \frac{e}{\varepsilon}\,(n_i-\ND), &&i=1,\dots,N, \label{poisson}\\
J_{i}-J_{i-1} &= -e\,\frac{dn_i}{dt}, &&i=1,\dots,N, \label{continuity}
\end{align}
where $F_i$ denotes the spatially-averaged electric field in the $i^{\mathrm{th}}$ period of the SL, $n_i$ is the 2D electron density in the $i^{\mathrm{th}}$ well, $e$ is the elementary charge, $\varepsilon$ is the dielectric constant,  $\ND=1.5\times 10^{11}\,$cm$^{-2}$ is the 2D doping density, and $J_i$ is the tunneling current density from the $i^{\mathrm{th}}$ to the \mbox{$(i+1)^{\mathrm{th}}$} well \footnote{The $0^{\mathrm{th}}$ and $(N+1)^{\mathrm{th}}$ wells refer to the two contact layers.}. Differentiating \eqref{poisson} with respect to time and inserting the result into \eqref{continuity} gives
\begin{equation} \label{Jtotd}
\varepsilon\,\frac{dF_i}{dt}+J_i = J(t) , \hspace*{5 mm} i=0,\dots,N,
\end{equation}
where the total current density $J(t)$ is the same for all periods.
The bias condition is
\begin{equation} \label{Vtot}
\sum_{i=0}^{N} F_i = \frac{V}{\ell}\,,
\end{equation}
where $V$ is the time-independent voltage bias across the entire SL (which is the control parameter for this system), and where $\ell$ is the distance between successive wells.

\begin{figure}[t]
\hspace{.3cm}\includegraphics{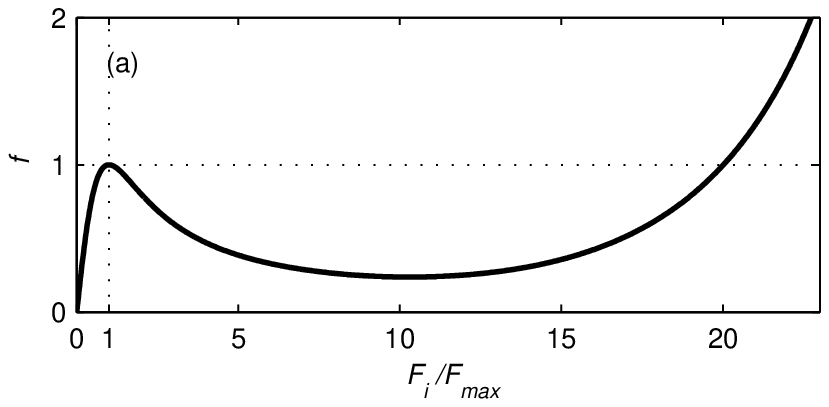}
\includegraphics{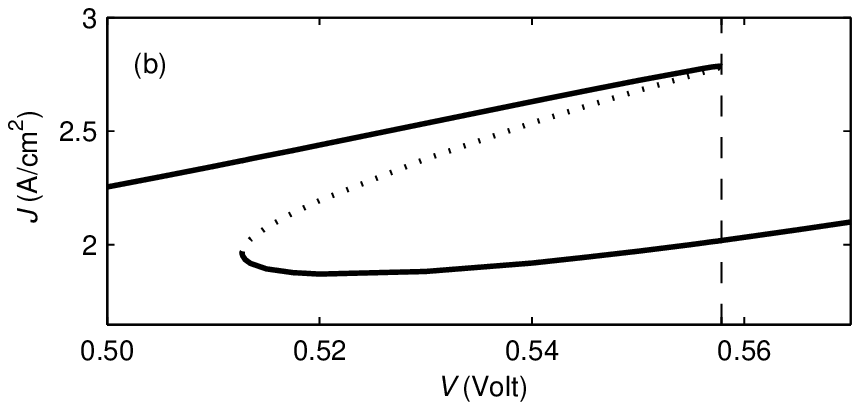}
\caption{
\label{fx JV fig}
(a) The $f(F_i/F_{\max})$ curve.
(b) The $J$-$V$ curve. Every point on the two solid branches corresponds to an attractor of our system, points on the dotted branch correspond to a saddle point. The vertical line at $\Vth$ marks the location of the saddle node bifurcation.}
\end{figure}

The current density $J_i$ is developed using a transfer Hamiltonian, and it can be written as a function of the local field and the neighboring charge densities \cite{XuTeitsworth}:
\begin{align}
  J_i\!
    &= \!\frac{e\,\vM f(F_i/F_{\max})}{\ell} \Big[n_i\!-\!c_1\ln\!\big\{1\!+e^{-c_2F_i}(e^{n_{i+1}/c_1}\!-\!1)\big\}\Big] \nonumber \\*
    &= J_i(F_i,n_i,n_{i+1})\label{SLmodel}
\end{align}
for $i=1,\dots,N-1$, where $c_1=1.68\cdot10^{10}$~$\mathrm{cm}^{\mathrm{-2}}$, $c_2=3.0145$\,cm/kV, $\vM=1.691$\,m/s, $F_{\max}=3.945$\,kV/cm, and where the function $f(F_i/F_{\max})$ is shown in Fig.~\ref{fx JV fig}(a)
\footnote{An exact expression for $f(\zeta)$ is derived in~\cite{XuTeitsworth}; here we use $f(\zeta)=2\zeta(1+\zeta^2)^{-1} + \exp\!\big(4\cdot10^{-6}\zeta^4\big) - 1$, a good approximation for the range of field values in this study.}.
%
%
%
%
%
%
For the present study, we use standard Ohmic boundary conditions
\begin{equation} \label{bc}
  J_0=g F_0
  \qquad\text{and}\qquad
  J_N=g F_N \frac{n_N}{\ND}\,,
\end{equation}
where $g=0.08$\,$(\Omega $m$)^{-1}$ denotes the effective contact \mbox{conductivity}.

Equations~\eqref{poisson} and \eqref{Jtotd}-\eqref{bc} constitute a complete dynamical system in $(F_0,\dots,F_N)$ describing electronic transport in the SL system and, for the parameters used here, are well-known to yield two stable points within a range of bias values $V$~\cite{XuTeitsworth}.
By \eqref{Jtotd}, \eqref{SLmodel}-\eqref{bc} and \eqref{poisson}, for each equilibrium solution we have $J_0=\dots=J_N=J=const$, and so the bistability can be illustrated by the current-voltage curve in Fig.\ \ref{fx JV fig}(b) which for every equilibrium point shows the associated value for~$J$.
The full curve consists of $N+1$ stable current branches, the first $N$ of which each ends in a saddle node bifurcation; Fig.\ \ref{fx JV fig}(b) shows only the $4^{\mathrm{th}}$ and $5^{\mathrm{th}}$ branches (the two solid curves), which are the focus of our study. Within the $n^{\mathrm{th}}$ current branch (counted from $V=0$), every stable state $(F_0,\dots,F_N)$ consists of two domains $(F_0,\dots,F_{N-n})$ and $(F_{N-n+1},\dots,F_N)$ within which the entries $F_i$ are roughly constant (low in the first and high in the second domain).
%

Let us now consider the effect of shot noise, which results from the discrete and random nature of the tunneling electrons \cite{BonillaSanchezSoler,Thesis_Xu}.
This noise induces rare transitions from the metastable state on branch $n$ to the one on branch $n+1$, which causes the overall current to jump and the high field domain to expand by one period.
Adding this noise to the tunneling current density, we modify~\eqref{Jtotd}~as
\begin{equation}
J(t)=\varepsilon\,\frac{dF_i}{dt}+J_i+J_i^{(r)}(t), \hspace*{5 mm} i=0,\dots,N.
\label{noiseJt}
\end{equation}
Using that the noise terms $J_i^{(r)}(t)$ are delta-correlated in time and space, one can solve \eqref{noiseJt} for $dF_i$ to write~\cite{Thesis_Xu}
\begin{align}
dF_i =\,& \frac{1}{\varepsilon} \Bigg[ \frac{1}{N+1}\sum_{j=0}^{N} J_j - J_i \Bigg]dt \nonumber \\*
&{}+ \frac{1}{\varepsilon}\Bigg[ \frac{1}{N+1}\sum_{j=0}^{N} B_j\,dW^j - B_i\,dW^i \Bigg]
\label{Emodel}
\end{align}
for $i=0,\dots,N$, where $B_i^2=\frac{e J_i}{a}$ (here $a$ denotes the cross-sectional area of the SL), and where the $W^i$ are independent Brownian motions.


\vspace{.2cm}
In order to apply gMAM to this model, we need to resolve the problem that the noise terms in the $(N+1)$-dimensional system \eqref{Emodel} are dependent (note that, because of \eqref{Vtot}, summing the second bracket in \eqref{Emodel} over all~$i$ gives $0$). We will therefore instead consider the $N$-dimensional transformed system for the variables~$n_i$,
\begin{align}
 \hspace{-8.5pt}dn_i &= \frac{\varepsilon}{e}\,(dF_i-dF_{i-1}) \nonumber\\
  &= \frac1e\big(J_{i-1}\!-\!J_{i}\big)\,dt
     +\frac1e\big(B_{i-1}\,dW^{i-1}-B_i\,dW^i\big) \label{ni system}\\
  &= \frac1e\big(J_{i-1}\!-\!J_{i}\big)\,dt
     +\!\frac1{\sqrt{ea}}\big(\sqrt{J_{i-1}}\,dW^{i-1}\!-\!\sqrt{J_i}\,dW^i\big)\hspace{3pt}
\nonumber
\end{align}
for $i=1,\dots,N$. To see that this is a closed system, observe that the linear system for the $F_i$ given by \eqref{poisson} and \eqref{Vtot} is solved by
\begin{equation*}
  F_i = \frac{V}{(N+1)\ell} + \frac e\varepsilon\Bigg[\Big(\frac N2 -i\Big)\ND + \sum_{j=1}^N\Big(\frac j{N+1}-\One_{j>i}\Big)n_j\Bigg]
\end{equation*}
for $i=0,\dots,N$, and so by \eqref{SLmodel}-\eqref{bc} we actually have $J_i=J_i(n_1,\dots,n_N)$ for $\forall i=0,\dots,N$.
Denoting $X:=(n_1,\dots,n_N)$, the system \eqref{ni system} is thus of the form
\begin{equation}\label{SDE}
  dX_t^\eta=b(X_t^\eta)\,dt+\sqrt\eta\,\sigma(X_t^\eta)\,dW_t,
\end{equation}
where $(W_t)_{t\geq0}$ is an $N$-dimensional Brownian motion, $\eta\!=\!(ea)^{-1}$ is a small parameter, and where the drift vector field $b(X)$ and the diffusion matrix $\sigma(X)$ are given~by
\begin{align}
  b(X) &=\frac1e\big(J_0-J_1,\dots,J_{N-1}-J_N\big), \label{drift field} \\[3pt]
  \sigma(X) &=
\begin{pmatrix}
 \sqrt{J_0} & -\sqrt{J_1}            & \hspace{19pt} 0         & \hspace{4pt}\cdots & \hspace{-7pt}0 \\[-2pt]
 0          & \hspace{7pt}\sqrt{J_1}   & -\sqrt{J_2}            & \hspace{4pt}\ddots & \hspace{-8pt}\vdots \\[-2pt]
 \vdots     &  \hspace{-14pt}\ddots                  & \hspace{-20pt}\ddots & \hspace{1pt}\ddots             &\hspace{-7pt}\raisebox{1pt}{0} \\[2pt]
 0          &  \hspace{-14pt}\cdots                  & \hspace{-35pt}0           & \hspace{-35pt}\sqrt{J_{N-1}} & \hspace{-7pt}-\sqrt{J_N}
\end{pmatrix}\!, \label{noise field}
\end{align}
with $J_i\!=\!J_i(X)$. As we will see below, this solves our problem of degenerate noise.


\vspace{.2cm}

One of the central results of large deviation theory is that the mean time for a transition from one metastable state $x_1$ to another, $x_2$, in a system of the form \eqref{SDE}, is given by $\langle T_\eta\rangle\approx e^{\hat U\hspace{-1pt}(x_1,x_2)/\eta}$ as $\eta\searrow0$, where $\hat U$ is the quasipotential \cite{WentzellFreidlin}.
%
%
It is well-known that in the gradient case (i.e., if $b=-\nabla U$ for some potential $U$, and if $\sigma$ is the identity matrix) we have $\hat U(x_1,x_2)=U(x_s)-U(x_1)$, where $x_s$ is the lowest-energy saddle point between the two basins of attraction; the maximum likelihood transition path is in this case the path of minimum energy, which can be numerically computed using the string method~\cite{String1}. However, the theory extends to the general case \eqref{SDE}, and the quasipotential~$\hat U$ is then given by
\begin{equation}\label{V geo}
  \hat U(x_1,x_2) = \inf_{\gamma\in\Gx}\int_\gamma\big(|b(x)|_x|dx|_x-\skpA{b(x)}{dx}\big),
\end{equation}
where $\Gx$ is the set of rectifiable (i.e., finite-length) curves leading from $x_1$ to $x_2$, and where $\skpA{u}{v}:=\skp{u}{A(x)^{-1}v}$ and $|u|_x:=\skpA{u}{u}^{1/2}$ for $\forall u,v\in\Rn$. Here we define $A(x):=\sigma(x)\sigma(x)^T$, which we require to be positive definite for every $x$.

Equation~\eqref{V geo} is a geometric reformulation \cite{CPAM_Hey1} of a more common formula for $\hat U$ that is based on time-parameterized paths \cite[Eqns.~(3.5) and (4.1)]{WentzellFreidlin}, and it has gotten some attention recently for its analytical and numerical advantages. In particular, the infimum in~\eqref{V geo} is typically achieved at some minimizing curve $\gamma^\star\in\Gx$ \cite{Existence1}, i.e., the \textit{maximum likelihood transition curve}, and so denoting the integral in \eqref{V geo} by $S(\gamma)$, we have
\begin{equation}
\hat U(x_1,x_2)=S(\gamma^\star).
\end{equation}

Furthermore, \eqref{V geo} is the basis of an efficient algorithm (the geometric Minimum Action Method, or gMAM) for computing~$\gamma^\star$.
The algorithm works by moving some initial curve $\gamma_0\in\Gx$ successively into the direction of steepest descent while keeping the end points fixed. If we write\vspace{-.1cm}
\[
  S(\gamma)=S(\varphi)=
\int_0^1\big(|b(\varphi)|_\varphi|\varphi'|_\varphi-\skp{b(\varphi)}{\varphi'}_\varphi\big)\,d\alpha
\vspace{-.1cm}
\]
for some parameterization $\varphi(\alpha)\colon[0,1]\to\Rn$ of $\gamma$ then gMAM numerically solves the PDE $\partial_\tau\varphi(\tau,\alpha)=-\delta S(\varphi)$, $(\tau,\alpha)\in[0,\infty)\times[0,1]$, subject to the boundary conditions $\varphi(\tau\!=\!0,\,\cdot\,)=\varphi_0$, $\varphi(\,\cdot\,,\alpha\!=\!0)=x_1$, $\varphi(\,\cdot\,,\alpha\!=\!1)=x_2$.
After preconditioning \cite{CPAM_Hey1}, this PDE can be written as
\begin{align}
 \partial_\tau\varphi
  =\,& \lambda^2\varphi'' - \lambda(\nabla b(\varphi)+C)\varphi' \nonumber \\
  &\hspace{1cm}+ A(\varphi) \big(\nabla b(\varphi)+\tfrac12 C\big)^T\theta
    + \lambda\lambda'\varphi', \label{PDE}
\end{align}
where $\varphi'$ and $\varphi''$ denote the $\alpha$-derivatives of~$\varphi$, \
$\lambda(\varphi,\varphi')$ $:=|b(\varphi)|_\varphi/|\varphi'|_\varphi$, \ $\theta(\varphi,\varphi'):=A(\varphi)^{-1}(\lambda\varphi'-b(\varphi))$, and
$C(\varphi,\varphi')$ is the $(N\times N)$-matrix whose $i^{\text{th}}$ column is $[\partial_{x_i}A]\big|_{x=\varphi}\,\theta$.
To stabilize the algorithm, a semi-implicit code is used for the $\varphi''$ term, and the discretization points are redistributed equidistantly along the curve after each iteration step. See \cite{CPAM_Hey1,PRL_Hey1} for the derivation of \eqref{PDE} in the additive case (i.e., if $A$ is the identity matrix) and for further details on the implementation of gMAM.


\vspace{.2cm}
To see that gMAM is applicable to our system \eqref{SDE}-\eqref{noise field}, first note that
\begin{equation*}
  A=\sigma\sigma^T=
\begin{pmatrix}
  J_0+J_1 & -J_1   &    \hspace{7pt}0      & \hspace{-17pt}\cdots\hspace{11pt} 0            \\[-2pt]
   -J_1 &  J_1+J_2 & \hspace{-0pt}-J_2     &   \hspace{-19pt}\raisebox{-3pt}{$\ddots$}\hspace{15pt}\raisebox{-2pt}{$\vdots$}    \\[0pt]
 0 & \hspace{-6pt}\raisebox{-10pt}{$\ddots$} & \hspace{-3pt}\raisebox{-8pt}{$\ddots$} & \hspace{-21pt}\raisebox{-3pt}{$\ddots$}\hspace{17pt}\raisebox{-2pt}{0} \\[-8pt]
  \vdots   & \hspace{-37pt}\ddots &  & -J_{N-1}    \\[2pt]
  0      & \hspace{-20pt}\cdots\hspace{10pt} 0 & \hspace{-7pt}-J_{N-1} & \hspace{0pt}J_{N-1}+J_N
\end{pmatrix};
\end{equation*}
it is then straightforward to show by induction on $N$ that $\det(A)=\sum_{k=0}^N\prod_{i=0,i\neq k}^{N}J_i$, which is positive since \mbox{$J_i>0$} for $\forall i=0,\dots,N$.

To apply gMAM to the SL system, we first fix $V=0.52$\,Volt (see Fig.~\ref{fx JV fig}(b)) and locate the two stable points $(n_1,\dots,n_N)$ of the transformed system with Newton's method. To find the saddle point in between, we make a novel application of the string method~\cite{String1}.
Briefly, the string method is a fast and easy algorithm for finding a curve that connects the two attractors and is everywhere parallel or antiparallel to the drift~$b$. While the original purpose of this algorithm was to find transition curves in \textit{gradient} systems (which have this geometric property),
in general systems one can still make use of the fact that by construction such a curve must contain an equilibrium point. One can therefore find the saddle point by (i) applying the string method, (ii) searching for an equilibrium point along the obtained curve, and (iii) increase its accuracy by using Newton's method.

Finally, we find the minimum action curve $\gamma^\star$ leading from the first attractor to the saddle point, and we compute the associated action $S(\gamma^\star)$.
(Note that, as a consequence of the non-gradient nature of this system, the curve $\gamma^\star$ is found to deviate significantly from the unstable manifold of the saddle point, see also \cite{ChanDykman}.)
By increasing $V$ and using continuity, we determine the three equilibrium points, and then $\gamma^\star$ and $S(\gamma^\star)$, all as functions of $V$. We identify $\Vth=0.5578288$\,Volt as the value~$V$ at which the attractor on the upper branch collides with the saddle point.



\begin{figure}[t]
\includegraphics{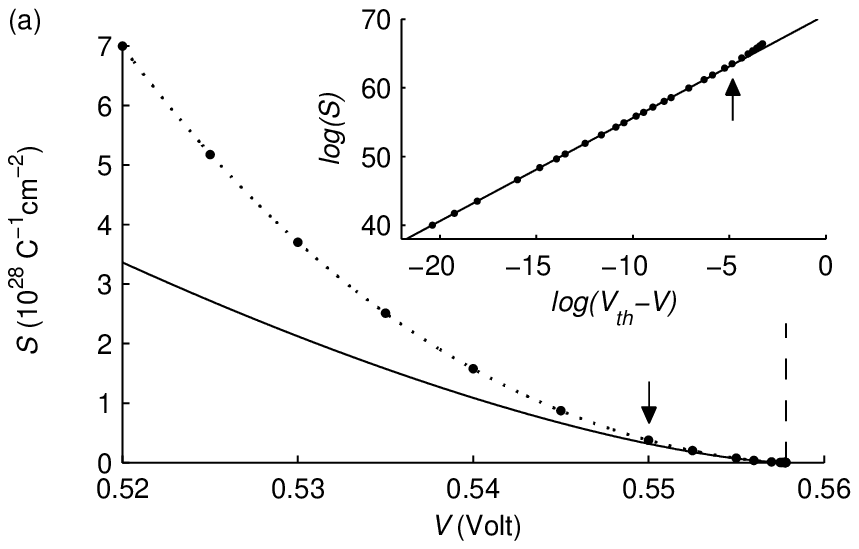}
\includegraphics{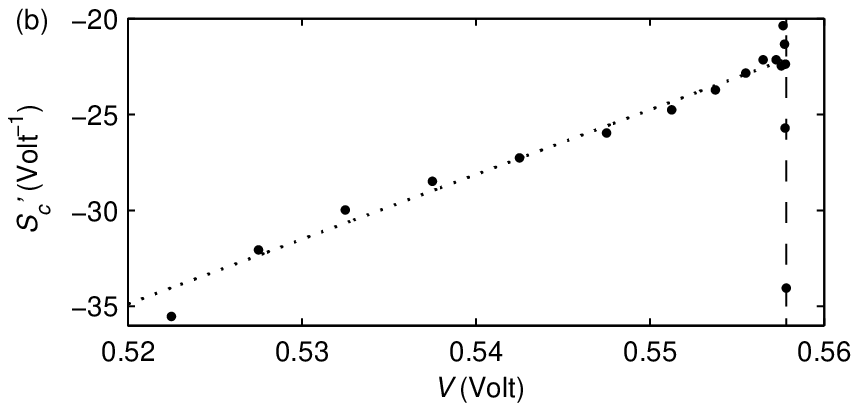}
\caption{\label{SV curve}(a) $S(V)$, (b) $S_c'(V)$. The solid circles are the values found by gMAM, solid and dotted curves are the approximations \eqref{fitted curve} and \eqref{fitted curve 2}, respectively. The vertical lines mark~$\Vth$.
}
\end{figure}

Fig.~\ref{SV curve}(a) shows plots of the action $S(V)=S(\gamma^\star(V))$. The values found by gMAM are shown as solid circles, the solid and dotted curves are the approximations found below. As the log-log plot in the inset illustrates, near the threshold~$\Vth$ we obtain the lowest-order approximation
\begin{equation} \label{fitted curve}
  S(V)\approx s_0 v^\beta ,\quad \text{where }v={|\Vth-V|}/{\Vth},
\end{equation}
$\beta=1.5\pm0.0001$ and $s_0=(1.905\pm0.0002)\cdot10^{30}\,\mathrm{C}^{-1}\mathrm{cm}^{-2}$.
%
The arrows mark the data point at which the slope of the log-log plot increases to $1.6$, and so we find that in this sense \eqref{fitted curve} is valid for $0.55\,\mathrm{Volt}\lesssim V\!\leq\Vth$.

To obtain higher-order terms, we consider $S_c(V):=S(V)/(s_0\,v^{3/2})$, which by construction fulfills $S_c(\Vth)=1$. The linear approximation of $S_c'(V)$ in Fig.\ \ref{SV curve}(b) leads to a quadratic approximation of $S_c(V)$, which yields
\begin{equation} \label{fitted curve 2}
  S(V)\approx s_0\,v^{3/2} + s_1\,v^{5/2} + s_2\,v^{7/2},
\end{equation}
where $s_1=2.35\cdot10^{31}\,\mathrm{C}^{-1}\mathrm{cm}^{-2}$ and $s_2=10^{32}\,\mathrm{C}^{-1}\mathrm{cm}^{-2}$.
This appears to be a good fit in the entire voltage range $0.52\,\mathrm{Volt}\leq V\!\leq\Vth$, see the dotted line in Fig.\ \ref{SV curve}(a).


\vspace{.2cm}
In this paper, we have generalized the gMAM algorithm for finding transition curves and their associated actions to SDE systems with multiplicative noise. We have used a new method for locating saddle points, and we have demonstrated the use of a geometric formulation \eqref{V geo} of the quasipotential.
\newline\indent
We have applied our techniques to the problem of stochastic current switching in semiconductor SLs, and to do so, we transformed the standard SL model into a form with non-degenerate noise matrix.
Our results for the SL show that in the zero-noise limit we have $\log\langle T\rangle\propto|\Vth-V|^{3/2}$ as $V\nearrow\Vth$, and we are able to quantify the range in which this scaling is valid. In addition, we have found a higher-order approximation that is valid also for larger values of $|\Vth-V|$.
\newline\indent
%
%
%
We thank Eric Vanden-Eijnden and Huidong Xu for helpful conversations.
This research was supported in part by NSF grant DMS-0616710 (M.H.\ and J.C.M.) and NSF grant DMR-0804232 (S.W.T.).


\bibliography{mybibliography}

\begin{thebibliography}{21}%
\makeatletter
\providecommand \@ifxundefined [1]{%
 \@ifx{#1\undefined}
}%
\providecommand \@ifnum [1]{%
 \ifnum #1\expandafter \@firstoftwo
 \else \expandafter \@secondoftwo
 \fi
}%
\providecommand \@ifx [1]{%
 \ifx #1\expandafter \@firstoftwo
 \else \expandafter \@secondoftwo
 \fi
}%
\providecommand \natexlab [1]{#1}%
\providecommand \enquote  [1]{``#1''}%
\providecommand \bibnamefont  [1]{#1}%
\providecommand \bibfnamefont [1]{#1}%
\providecommand \citenamefont [1]{#1}%
\providecommand \href@noop [0]{\@secondoftwo}%
\providecommand \href [0]{\begingroup \@sanitize@url \@href}%
\providecommand \@href[1]{\@@startlink{#1}\@@href}%
\providecommand \@@href[1]{\endgroup#1\@@endlink}%
\providecommand \@sanitize@url [0]{\catcode `\\12\catcode `\$12\catcode
  `\&12\catcode `\#12\catcode `\^12\catcode `\_12\catcode `\%12\relax}%
\providecommand \@@startlink[1]{}%
\providecommand \@@endlink[0]{}%
\providecommand \url  [0]{\begingroup\@sanitize@url \@url }%
\providecommand \@url [1]{\endgroup\@href {#1}{\urlprefix }}%
\providecommand \urlprefix  [0]{URL }%
\providecommand \Eprint [0]{\href }%
\@ifxundefined \urlstyle {%
  \providecommand \doi  [0]{\begingroup \@sanitize@url \@doi}%
  \providecommand \@doi [1]{\endgroup \@@startlink {\doibase
  #1}doi:\discretionary {}{}{}#1\@@endlink }%
}{%
  \providecommand \doi  [0]{doi:\discretionary{}{}{}\begingroup
  \urlstyle{rm}\Url }%
}%
\providecommand \doibase [0]{http://dx.doi.org/}%
\providecommand \Doi [0]{\begingroup \@sanitize@url \@Doi }%
\providecommand \@Doi  [1]{\endgroup\@@startlink{\doibase#1}\@@Doi}%
\providecommand \@@Doi [1]{#1\@@endlink}%
\providecommand \selectlanguage [0]{\@gobble}%
\providecommand \bibinfo  [0]{\@secondoftwo}%
\providecommand \bibfield  [0]{\@secondoftwo}%
\providecommand \translation [1]{[#1]}%
\providecommand \BibitemOpen [0]{}%
\providecommand \bibitemStop [0]{}%
\providecommand \bibitemNoStop [0]{.\EOS\space}%
\providecommand \EOS [0]{\spacefactor3000\relax}%
\providecommand \BibitemShut  [1]{\csname bibitem#1\endcsname}%
\bibitem [{\citenamefont {den Otter}(2000)}]{denOtter}%
  \BibitemOpen
  \bibfield  {author} {\bibinfo {author} {\bibfnamefont {W.~K.}\ \bibnamefont
  {den Otter}},\ }\href@noop {} {\bibfield  {journal} {\bibinfo  {journal} {J.
  Chem. Phys.},\ }\textbf {\bibinfo {volume} {112}},\ \bibinfo {pages} {7283}
  (\bibinfo {year} {2000})}\BibitemShut {NoStop}%
\bibitem [{\citenamefont {Aldridge}\ and\ \citenamefont
  {Cleland}(2005)}]{AldridgeCleland}%
  \BibitemOpen
  \bibfield  {author} {\bibinfo {author} {\bibfnamefont {J.~S.}\ \bibnamefont
  {Aldridge}}\ and\ \bibinfo {author} {\bibfnamefont {A.~N.}\ \bibnamefont
  {Cleland}},\ }\href@noop {} {\bibfield  {journal} {\bibinfo  {journal} {Phys.
  Rev. Lett.},\ }\textbf {\bibinfo {volume} {94}},\ \bibinfo {pages} {156403}
  (\bibinfo {year} {2005})}\BibitemShut {NoStop}%
\bibitem [{\citenamefont {Chan}\ and\ \citenamefont
  {Stambaugh}(2007)}]{ChanStambaugh}%
  \BibitemOpen
  \bibfield  {author} {\bibinfo {author} {\bibfnamefont {H.~B.}\ \bibnamefont
  {Chan}}\ and\ \bibinfo {author} {\bibfnamefont {C.}~\bibnamefont
  {Stambaugh}},\ }\href@noop {} {\bibfield  {journal} {\bibinfo  {journal}
  {Phys. Rev. Lett.},\ }\textbf {\bibinfo {volume} {99}},\ \bibinfo {pages}
  {060601} (\bibinfo {year} {2007})}\BibitemShut {NoStop}%
\bibitem [{\citenamefont {Chan}\ \emph {et~al.}(2008)\citenamefont {Chan},
  \citenamefont {Dykman},\ and\ \citenamefont {Stambaugh}}]{ChanDykman}%
  \BibitemOpen
  \bibfield  {author} {\bibinfo {author} {\bibfnamefont {H.~B.}\ \bibnamefont
  {Chan}}, \bibinfo {author} {\bibfnamefont {M.~I.}\ \bibnamefont {Dykman}}, \
  and\ \bibinfo {author} {\bibfnamefont {C.}~\bibnamefont {Stambaugh}},\
  }\href@noop {} {\bibfield  {journal} {\bibinfo  {journal} {Phys. Rev.
  Lett.},\ }\textbf {\bibinfo {volume} {100}},\ \bibinfo {pages} {130602}
  (\bibinfo {year} {2008})}\BibitemShut {NoStop}%
\bibitem [{\citenamefont {Cui}\ \emph {et~al.}(2010)\citenamefont {Cui} \emph
  {et~al.}}]{Cui}%
  \BibitemOpen
  \bibfield  {author} {\bibinfo {author} {\bibfnamefont {Y.-T.}\ \bibnamefont
  {Cui}} \emph {et~al.},\ }\href@noop {} {\bibfield  {journal} {\bibinfo
  {journal} {Phys. Rev. Lett.},\ }\textbf {\bibinfo {volume} {104}},\ \bibinfo
  {pages} {097201} (\bibinfo {year} {2010})}\BibitemShut {NoStop}%
\bibitem [{\citenamefont {Allen}\ \emph {et~al.}(2005)\citenamefont {Allen},
  \citenamefont {Warren},\ and\ \citenamefont {ten Wolde}}]{AllenWarren}%
  \BibitemOpen
  \bibfield  {author} {\bibinfo {author} {\bibfnamefont {R.~J.}\ \bibnamefont
  {Allen}}, \bibinfo {author} {\bibfnamefont {P.~B.}\ \bibnamefont {Warren}}, \
  and\ \bibinfo {author} {\bibfnamefont {P.~R.}\ \bibnamefont {ten Wolde}},\
  }\href@noop {} {\bibfield  {journal} {\bibinfo  {journal} {Phys. Rev.
  Lett.},\ }\textbf {\bibinfo {volume} {94}},\ \bibinfo {pages} {018104}
  (\bibinfo {year} {2005})}\BibitemShut {NoStop}%
\bibitem [{\citenamefont {E}\ \emph {et~al.}(2002)\citenamefont {E},
  \citenamefont {Ren},\ and\ \citenamefont {Vanden-Eijnden}}]{String1}%
  \BibitemOpen
  \bibfield  {author} {\bibinfo {author} {\bibfnamefont {W.}~\bibnamefont {E}},
  \bibinfo {author} {\bibfnamefont {W.}~\bibnamefont {Ren}}, \ and\ \bibinfo
  {author} {\bibfnamefont {E.}~\bibnamefont {Vanden-Eijnden}},\ }\href@noop {}
  {\bibfield  {journal} {\bibinfo  {journal} {Physical Review B},\ }\textbf
  {\bibinfo {volume} {66}},\ \bibinfo {pages} {052301} (\bibinfo {year}
  {2002})}\BibitemShut {NoStop}%
\bibitem [{\citenamefont {E}\ \emph {et~al.}(2004)\citenamefont {E},
  \citenamefont {Ren},\ and\ \citenamefont {Vanden-Eijnden}}]{MAM}%
  \BibitemOpen
  \bibfield  {author} {\bibinfo {author} {\bibfnamefont {W.}~\bibnamefont {E}},
  \bibinfo {author} {\bibfnamefont {W.}~\bibnamefont {Ren}}, \ and\ \bibinfo
  {author} {\bibfnamefont {E.}~\bibnamefont {Vanden-Eijnden}},\ }\href@noop {}
  {\bibfield  {journal} {\bibinfo  {journal} {Communications in Pure and
  Applied Mathematics},\ }\textbf {\bibinfo {volume} {57}},\ \bibinfo {pages}
  {637} (\bibinfo {year} {2004})}\BibitemShut {NoStop}%
\bibitem [{\citenamefont {Heymann}\ and\ \citenamefont
  {Vanden-Eijnden}(2008){\natexlab{a}}}]{CPAM_Hey1}%
  \BibitemOpen
  \bibfield  {author} {\bibinfo {author} {\bibfnamefont {M.}~\bibnamefont
  {Heymann}}\ and\ \bibinfo {author} {\bibfnamefont {E.}~\bibnamefont
  {Vanden-Eijnden}},\ }\href@noop {} {\bibfield  {journal} {\bibinfo  {journal}
  {Communications in Pure and Applied Mathematics},\ }\textbf {\bibinfo
  {volume} {61}},\ \bibinfo {pages} {1052} (\bibinfo {year}
  {2008}{\natexlab{a}})}\BibitemShut {NoStop}%
\bibitem [{\citenamefont {Heymann}\ and\ \citenamefont
  {Vanden-Eijnden}(2008){\natexlab{b}}}]{PRL_Hey1}%
  \BibitemOpen
  \bibfield  {author} {\bibinfo {author} {\bibfnamefont {M.}~\bibnamefont
  {Heymann}}\ and\ \bibinfo {author} {\bibfnamefont {E.}~\bibnamefont
  {Vanden-Eijnden}},\ }\href@noop {} {\bibfield  {journal} {\bibinfo  {journal}
  {Phys. Rev. Lett.},\ }\textbf {\bibinfo {volume} {100}},\ \bibinfo {pages}
  {140601} (\bibinfo {year} {2008}{\natexlab{b}})}\BibitemShut {NoStop}%
\bibitem [{\citenamefont {Vanden-Eijnden}\ and\ \citenamefont
  {Heymann}(2008)}]{JCP_Hey1}%
  \BibitemOpen
  \bibfield  {author} {\bibinfo {author} {\bibfnamefont {E.}~\bibnamefont
  {Vanden-Eijnden}}\ and\ \bibinfo {author} {\bibfnamefont {M.}~\bibnamefont
  {Heymann}},\ }\href@noop {} {\bibfield  {journal} {\bibinfo  {journal}
  {Journal of Chemical Physics},\ }\textbf {\bibinfo {volume} {128}},\ \bibinfo
  {pages} {061103} (\bibinfo {year} {2008})}\BibitemShut {NoStop}%
\bibitem [{\citenamefont {Freidlin}\ and\ \citenamefont
  {Wentzell}(1998)}]{WentzellFreidlin}%
  \BibitemOpen
  \bibfield  {author} {\bibinfo {author} {\bibfnamefont {M.~I.}\ \bibnamefont
  {Freidlin}}\ and\ \bibinfo {author} {\bibfnamefont {A.~D.}\ \bibnamefont
  {Wentzell}},\ }\href@noop {} {\emph {\bibinfo {title}
  {\href{http://www.amazon.com/Random-Perturbations-Dynamical-Systems-Freidlin%
/dp/3540983627}{Random Perturbations}
  \href{http://www.amazon.com/Random-Perturbations-Dynamical-Systems-Freidlin/%
dp/3540983627}{of Dynamical Systems}}}},\ \bibinfo {edition} {2nd}\ ed.\
  (\bibinfo  {publisher} {Springer, NY},\ \bibinfo {year} {1998})\BibitemShut
  {NoStop}%
\bibitem [{\citenamefont {Heymann}(2010)}]{Existence1}%
  \BibitemOpen
  \bibfield  {author} {\bibinfo {author} {\bibfnamefont {M.}~\bibnamefont
  {Heymann}},\ }\href@noop {} {\bibfield  {journal} {\bibinfo  {journal}
  {arXiv:1004.4873v1}} (\bibinfo {year} {2010})}\BibitemShut {NoStop}%
\bibitem [{\citenamefont {Bonilla}\ and\ \citenamefont
  {Teitsworth}(2010)}]{BonillaTeitsworth}%
  \BibitemOpen
  \bibfield  {author} {\bibinfo {author} {\bibfnamefont {L.~L.}\ \bibnamefont
  {Bonilla}}\ and\ \bibinfo {author} {\bibfnamefont {S.~W.}\ \bibnamefont
  {Teitsworth}},\ }\href@noop {} {\emph {\bibinfo {title}
  {\href{http://www.wiley-vch.de/publish/en/books/ISBN978-3-527-40695-1}{Nonli%
near Wave}
  \href{http://www.wiley-vch.de/publish/en/books/ISBN978-3-527-40695-1}{Methods
  for Charge Transport}}}},\ \bibinfo {edition} {1st}\ ed.\ (\bibinfo
  {publisher} {Wiley-VCH, Berlin},\ \bibinfo {year} {2010})\BibitemShut
  {NoStop}%
\bibitem [{\citenamefont {Bonilla}\ \emph {et~al.}(2002)\citenamefont
  {Bonilla}, \citenamefont {S\'anchez},\ and\ \citenamefont
  {Soler}}]{BonillaSanchezSoler}%
  \BibitemOpen
  \bibfield  {author} {\bibinfo {author} {\bibfnamefont {L.~L.}\ \bibnamefont
  {Bonilla}}, \bibinfo {author} {\bibfnamefont {O.}~\bibnamefont {S\'anchez}},
  \ and\ \bibinfo {author} {\bibfnamefont {J.}~\bibnamefont {Soler}},\
  }\href@noop {} {\bibfield  {journal} {\bibinfo  {journal} {Physical Review
  B},\ }\textbf {\bibinfo {volume} {65}},\ \bibinfo {pages} {195308} (\bibinfo
  {year} {2002})}\BibitemShut {NoStop}%
\bibitem [{\citenamefont {Rogozia}\ \emph {et~al.}(2001)\citenamefont {Rogozia}
  \emph {et~al.}}]{Rogozia}%
  \BibitemOpen
  \bibfield  {author} {\bibinfo {author} {\bibfnamefont {M.}~\bibnamefont
  {Rogozia}} \emph {et~al.},\ }\href@noop {} {\bibfield  {journal} {\bibinfo
  {journal} {Phys. Rev. B},\ }\textbf {\bibinfo {volume} {64}},\ \bibinfo
  {pages} {041308} (\bibinfo {year} {2001})}\BibitemShut {NoStop}%
\bibitem [{\citenamefont {Tretiakov}\ \emph {et~al.}(2003)\citenamefont
  {Tretiakov}, \citenamefont {Gramespacher},\ and\ \citenamefont
  {Matveev}}]{Tretiakov}%
  \BibitemOpen
  \bibfield  {author} {\bibinfo {author} {\bibfnamefont {O.~A.}\ \bibnamefont
  {Tretiakov}}, \bibinfo {author} {\bibfnamefont {T.}~\bibnamefont
  {Gramespacher}}, \ and\ \bibinfo {author} {\bibfnamefont {K.~A.}\
  \bibnamefont {Matveev}},\ }\href@noop {} {\bibfield  {journal} {\bibinfo
  {journal} {Phys. Rev. B},\ }\textbf {\bibinfo {volume} {67}},\ \bibinfo
  {pages} {073303} (\bibinfo {year} {2003})}\BibitemShut {NoStop}%
\bibitem [{\citenamefont {Xu}\ and\ \citenamefont
  {Teitsworth}(2007)}]{XuTeitsworth}%
  \BibitemOpen
  \bibfield  {author} {\bibinfo {author} {\bibfnamefont {H.}~\bibnamefont
  {Xu}}\ and\ \bibinfo {author} {\bibfnamefont {S.~W.}\ \bibnamefont
  {Teitsworth}},\ }\href@noop {} {\bibfield  {journal} {\bibinfo  {journal}
  {Physical Review B},\ }\textbf {\bibinfo {volume} {76}},\ \bibinfo {pages}
  {235302} (\bibinfo {year} {2007})}\BibitemShut {NoStop}%
\bibitem [{Note1()}]{Note1}%
  \BibitemOpen
  \bibinfo {note} {The $0^{\protect \mathrm {th}}$ and $(N+1)^{\protect \mathrm
  {th}}$ wells refer to the two contact layers.}\BibitemShut {Stop}%
\bibitem [{Note2()}]{Note2}%
  \BibitemOpen
  \bibinfo {note} {An exact expression for $f(\zeta )$ is derived in~\cite
  {XuTeitsworth}; here we use $f(\zeta )=2\zeta (1+\zeta ^2)^{-1} + \protect
  \qopname \relax o{exp}\protect \tmspace -\thinmuskip {.1667em}{\setbox \z@
  \hbox {\frozen@everymath \@emptytoks \mathsurround \z@ $\nulldelimiterspace
  \z@ \left (\vcenter to\@ne \big@size {}\right .$}\box \z@ }4\cdot
  10^{-6}\zeta ^4{\setbox \z@ \hbox {\frozen@everymath \@emptytoks
  \mathsurround \z@ $\nulldelimiterspace \z@ \left )\vcenter to\@ne \big@size
  {}\right .$}\box \z@ } - 1$, a good approximation for the range of field
  values in this study.}\BibitemShut {Stop}%
\bibitem [{\citenamefont {Xu}(2010)}]{Thesis_Xu}%
  \BibitemOpen
  \bibfield  {author} {\bibinfo {author} {\bibfnamefont {H.}~\bibnamefont
  {Xu}},\ }\href@noop {} {Ph.D. thesis},\ \bibinfo  {school} {Duke University}
  (\bibinfo {year} {2010})\BibitemShut {NoStop}%
\end{thebibliography}%

\end{document}